\documentclass{article}
\usepackage{amssymb}
\usepackage{graphicx}
\usepackage{amsmath}

\newtheorem{theorem}{Theorem}[section]

\newtheorem{conjecture}[theorem]{Conjecture}

\newtheorem{definition}[theorem]{Definition}

\newcommand{\SL}{\mathrm{SL}}

\newcommand{\mr}[1]{\mathcal{#1}}

\newcommand{\Z}{\mathbb{Z}}

\newcommand{\F}{\mathbb{F}}

\newcommand{\KC}{\mr{K}}

\newcommand{\KaC}{Kazhdan constant}
\newcommand{\Sym}{\mathrm{Sym}}
\newcommand{\Alt}{\mathrm{Alt}}

\begin{document}
\title{Symmetric Groups and Expanders}
\author{M. Kassabov } \date{}
\author{Martin Kassabov}

\maketitle
{
\renewcommand{\thefootnote}{}
}
\begin{abstract}
We construct an explicit generating sets $F_n$ and $\tilde F_n$ of the 
alternating and the symmetric groups,
which make the Cayley graphs $\mr{C}(\Alt(n), F_n)$ and $\mr{C}(\Sym(n), \tilde F_n)$
a family of bounded degree expanders for all sufficiently large $n$.
These expanders have many applications
in the theory of random walks on groups and other areas of mathematics.
\end{abstract}

\renewcommand{\thetheorem}{\arabic{theorem}}

A finite graph $\Gamma$ is called an $\epsilon$-expander for some
$\epsilon \in (0,1)$, if for any subset $A \subseteq \Gamma$ of size
at most $|\Gamma|/2$ we have $|\partial (A)| >\epsilon |A|$
(where $\partial(A)$ is the set of vertices of $\Gamma \backslash A$ of edge
distance 1 to $A$). The largest such $\epsilon$ is called the
expanding constant of $\Gamma$. Constructing families of $\epsilon$-expanders
with 
bounded valency
is an important practical problem in computer science,
because such graphs have many nice properties --- for example they
have a logarithmic diameter. For an excellent
introduction to the subject we refer the reader to the
book~\cite{expanderbook} by A. Lubotzky.

Using counting arguments it can be shown that almost
any $5$ regular graph is $1/5$-expander. However constructing
an explicit examples of families expander graphs is
a difficult problem.

\medskip

The first explicit construction of a family of expanders was
done by G. Margulis in~\cite{Mar}, using Kazhdan property \emph{T}
of $\SL_3(\Z)$. Currently there are several different construction
of expanders. With the exception of a few recent ones based on
the zig-zag products of graphs
(see~\cite{ALW,RVW,RSW}), all constructions are based groups theory and use
some variant of property \emph{T} (property $\tau$, Selberg
property etc.).

Kazhdan Property \emph{T} is not very interesting for a given finite group $G$
(all finite groups have property \emph{T}),
but the related Kazhdan constant with respect to some generating $F$
set is. Given an infinite collection of finite groups $G_i$, it is a
challenge to prove the existence of uniform Kazhdan constants
with respect to properly chosen generating sets. This problem is
related to construction a family of expanders using the Cayley graphs
of the groups $G_i$.

The original definition of property \emph{T} uses the Fell
topology of the unitary dual, see~\cite{kazhdan}.
Here we will use an equivalent
definition (only for discrete groups) which also addresses the
notion of the \KaC s.

\begin{definition}
\label{expander}
Let $G$ be a discrete group generated by a finite set $S$.
Then $G$ has the Kazhdan property \emph{T} if there exists
$\epsilon  >0$
such that for every unitary
representation $\rho: \ G \rightarrow U(\cal{H})$ on a
Hilbert space $\mr{H}$ without $G$ invariant vectors and
every  vector $v \not = 0$ there
exists some $s \in S$ such that $||\rho(s)v-v||> \epsilon ||v||$.
The largest $\epsilon$ with this property is called the \emph{\KaC}
for $G$ with respect to $S$ and is denoted by $\KC(G;S)$.
\end{definition}
For a 
$G$ the property \emph{T} is independent on the
choice of the generating set $S$, however the Kazhdan constant
depends also on the generating set.

\medskip

The following connection between property \emph{T} and expander graphs
is well known:

\begin{theorem}%
[\cite{expanderbook}, Theorem 4.3.2]
\label{CGexpanders}
Let $G$ be a discrete group having property \emph{T}, and let
$S$ be a finite generating set of $G$.
Then there exists an
$\epsilon=\epsilon(S) >0$ such that
the Cayley graphs $\mr{C}(G_i,S_i)$ (and all their quotients) of the finite
images of $G_i$ of $G$ (with respect to the images $S_i$ of $S$)
form a family of $\epsilon$-expanders. The largest $\epsilon_0(S)$
with this property is related to the \KaC\ $\KC(G,S)$,
in particular we have
$\epsilon_0(S) \geq \KC(G;S)^2/4$.
\end{theorem}

\medskip

Using this approach and property \emph{T} of $\SL_n(\Z)$
one can make the Cayley graphs of $\SL_n(\F_p)$ for fixed $n\geq 3$ a family
of expanders
.\footnote{
The Cayley graphs of $\SL_2(\F_p)$ can also be made expanders although
the group $\SL_2(\Z)$ does not have property \emph{T} or
even $\tau$, see~\cite{LPS}.}
Until recently the only way to
prove Kazhdan property \emph{T} was via representation theory of
high rank Lie groups. These methods are not
quantitative and does not lead to estimates for the
\KaC s and the corresponding expanding constants of
the resulting Cayley graphs.

A breakthrough in this direction was done in~\cite{YSh} by Y.~Shalom, who
used the bounded generation of the $\SL_n(\Z)$ and
M.~Burger's estimate (see~\cite{Bur})
of the relative \KaC\ of $\SL_2(\Z)\ltimes \Z^2$
to obtain an estimates for the \KaC\ of $\SL_n(\Z)$. His methods were refined
in~\cite{K} to give the exact asymptotic of the \KaC\ of $\SL_n(\Z)$
with respect to the set of all elementary matrices.
This yields an asymptotically exact estimate for the expansion
constant of the form $O(1/n)$ of the Cayley graphs of
$\SL_n(\F_p)$ with respect the set of all elementary
matrices.

It is interesting to note that if the rank of the matrices
increases then the resulting Cayley graphs do not form an
expander family, even though that the degree of these graphs
goes to infinity.

\medskip

Using relative property \emph{T} of the pair
$\SL_2(R)\ltimes R^2, R^2$ for finitely generated noncommutative rings
$R$, the Cayley graphs of $\SL_n(\F_q)$ for any prime power $q$
and infinitely many $n$-es can be
made expanders by choosing a suitable generating
set, see~\cite{KSL3k}. An important building block
in this construction is that the group $\SL_n(\F_q)$
can be written as a product of $20$ abelian subgroups and this number is
independent on $n$ and $q$.

However, the symmetric/alternating groups can not be written
as a product of fixed number of abelian subgroups,
because the size of the $\Sym(n)$ or $\Alt(n)$ is
approximately $n^{n}$ and every abelian subgroup
has at most $2^n$ elements, see~\cite{A}.
This suggests that $\Alt(n)$ are further from the abelian groups
than all others finite simple groups, and therefore
they should have more expanding properties.
Unfortunately, this also
significantly complicates the construction of expanders
based on the alternating groups.

\medskip

Using the classification of the finite simple groups A.~Lubotzky suggested~\cite{Lp}
that the results from~\cite{KSL3k} could be generalized to:
\begin{conjecture}
\label{FSGE}
There exists constants $L>0$ and $\epsilon>0$ such that for any
non-abelian finite simple group $G$, there exits a
generating set $F$ such that $|F|\leq L$ and
the Cayley graphs $\mr{C}(G;F)$ from a family of
$\epsilon$-expanders. Equivalently, we have that
$\KC(G;F)>\epsilon$ (with a different $\epsilon$).
\end{conjecture}
This conjecture is supported by
several results --- it is known (see~\cite{BLK} and~\cite{KR})
that for any non-abelian finite simple
group there exist a $4$ element generating set such that the
diameter of the corresponding Cayley graph is logarithmic in
the size of the group. Theorem~\ref{main1} together
with the results form~\cite{KSL3k} can be
view as a major step toward proving Conjecture~\ref{FSGE}, see~\cite{KNFS}.

\medskip

If one allows unbounded generating sets there are only
a few partial results known. A classical result of N.~Alon and Y.~Roichman
(and its improvements in~\cite{LR} and~\cite{LS}) says that
any group is an $\epsilon$-expander%
\footnote{There are two different definitions of $\epsilon$-expanders,
which are equivalent for graphs of bounded degree but not in general.
A `weak' one corresponding to Definition~\ref{expander},
and a `strong' one using the spectral gap of the Laplacian.
The random Cayley graphs obtained from Theorem~\ref{AlRo} are
expanders in both definitions.}
with respect to a random large generating set:
\begin{theorem}[\cite{AR}]
\label{AlRo}
For any $\epsilon>0$ there exists $c(\epsilon)>0$ such that
the Cayley graph any finite group $G$ is an $\epsilon$-expander
with respect to a random generating set of size
$c(\epsilon) \log |G|$.
\end{theorem}
The bound $c(\epsilon) \log |G|$ is the optimal one in the
class of all finite groups, because large the abelian groups
are not expanders with respect to small generating sets.
It is believed that if the group $G$ is far from being abelian
the bound $c\log |G|$ can be improved.
Note that Theorem~\ref{AlRo}
implies a version of Conjecture~\ref{FSGE}, where the
size of the generating set is allowed to increase.

Theorem~\ref{AlRo} also gives that the Cayley graphs of the
alternating groups $\Alt(n)$ are expanders with respect to a
random generating set of size $n \log n$.
This was the best known result in this direction and even there
were no known explicit sets of size less than
$n^{\sqrt{n}}$ which make the Cayley graphs expanders. In view of the
main results in this paper it is very interesting to understand
the expanding properties of the Cayley graphs of $\Alt(n)$  and
other finite simple groups with respect
to a random generating set of a small (even bounded) size.

\medskip

One of the main results in this paper gives
answers affirmatively an old question, which have been asked
several times in the literature, see~\cite{BHKLS,expanderbook,lu}:
\begin{theorem}
\label{main}
There exist constants $L>0$, $\epsilon >0$ and an infinite
sequence $n_s$ with the property:
There exists a constructible generating set
$F_{n_s}$ of size at most $L$ of the alternating group $\Alt(n_s)$ such that
the Cayley graphs $\mr{C}(\Alt(n_s);F_{n_s})$ form a family of
$\epsilon$-expanders.
\end{theorem}

From the proof of Theorem~\ref{main}, it can be seen that
the sequence $n_i$ does not
grow too fast, which leads to the following generalization:

\begin{theorem}
\label{main1}
There exist constants $L>0$ and $\epsilon >0$, with the
property: for any $n$ there exists
a constructible generating sets $F_{n}$ and $\tilde F_{n}$
of the alternating group $\Alt(n)$
and the symmetric group $\Sym(n)$ respectively such that
the Cayley graphs $\mr{C}(\Alt(n);F_{n})$ and
$\mr{C}(\Sym(n);\tilde F_{n})$ form a family of
$\epsilon$-expanders and all generating sets $F_n$ and
$\tilde F_n$ have at most $L$ elements.
\end{theorem}

\bigskip

Theorem~\ref{main1} has applications interesting applications:
It provides one of the few constructions of an expander family
of Cayley graphs $\mr{C}(G_i,S_i)$ such that the groups $G_i$
are not quotients of some infinite group having a variant
of Kazhdan property \emph{T}. It also provides a supporting evidence
for the conjecture that the automorphism groups of the free
groups have property $\tau$, see~\cite{LP}.

\medskip

Theorem~\ref{main1} implies that the expanding constant of
$\Alt(n)$ with respect to the set $F_n^{10^{10}}$ is large enough.%
\footnote{More precisely we have that the spectral gap of the Laplacian of
the Cayley graph is very close to $1$.}
The size of the set $F_n^{10^{10}}$ is independent on $n$, and
if $n$ is sufficiently big then $|F_n^{10^{10}}| < n^{1/30}/10^{-30}$.
The last inequality allows us to use the expander
$\mr{C}(\Alt(n);F_n^{10^{10}})$ as a
`seed' graph for E.~Rozemann, A.~Shalev and A.~Widgerson
recursive construction of expanders, see~\cite{RSW}.
This construction produces a family of expander graphs
based on Cayley graphs of the automorphism group of
large $n$-regular rooted tree of depth $k$. A slight modification of
this construction gives an other recursive expander family based on
$\Alt(n^k)$ for fixed large $n$ and different $k$-es.

\medskip

Theorem~\ref{main1} gives that the Cayley graphs $\mr{C}(\Alt(n); F_n)$ and
$\mr{C}(\Sym(n); \tilde F_n)$ have many expanding properties which
imply that the random walks on $\Alt(n)$ and $\Sym(n)$,
generated by $F_n$ and $\tilde F_n$
respectively, have mixing time approximately 
$\log |\Alt(n) | = n \log n$ steps. This leads to a natural and
fast algorithm for generating pseudo-random permutations.

\bigskip

\noindent
\textbf{Proof of Theorem~\ref{main}:}
We will think that the alternating group $\Alt(N)$ acts
on a set of $N$ points which are arranged into $d$ dimensional
cube of size $K$. Let $\rho$ be a representation of $\Alt(N)$ with
almost invariant vector $v$ with respect to some generating set $F$
(to be chosen later). We will use two different arguments to show that
there is an invariant vector.

First, we will break the representation $\rho$ into
two components --- one corresponding
to partitions $\lambda$ with $\lambda_1 < N - h$ and
second one containing all other partitions.
The decomposition of the regular representation of $\Alt(N)$ into two
components depending on the first part of the partition $\lambda$
comes from~\cite{Ro1}. In this paper, Y.~Roichman uses similar
argument to show that the Cayley
graphs of the symmetric/alternating group with respect to a conjugacy
class with a large number of
non-fixed points have certain expanding properties.

We will show that the projection of the vector $v$ in the first representation
is small provided that $h \gg K$. Also the projection of $v$
in the second one is close to an
invariant vector if $h \ll N^{1/4}$.%
\footnote{We believe that the argument also works even in the case $h \ll N^{1/2}$,
however we are unable to prove it.
If such generalization is true, it will allow us to use $d=4$,
which will improve the estimates
for the \KaC s by a factor of $10$.}
In order to satisfy these restriction we need that
$N \gg K^4$, i.e., $d>4$. In order to simplify the argument a little, we also
require that $d$ is even, which justifies our choice of $d=6$ and $N=K^6$.
Also we need that $K \ll h \ll K^{3/2}$, therefore using $h=K^{5/4}$
we will define $\mr{H}_1$ to be the sub-representation of
$\mr{H}$ corresponding to all partitions
with $\lambda_1 < N - K^{5/4}$.
As an additional assumption, we require that $K+1$ is a
power of some prime number and we will use
$K=2^{3s}-1$ for a significantly large $s$.%
\footnote{It seems that $s>50$ suffices, however the estimate depends on
the constants involved in the character estimates from~\cite{Ro},
which are not in the literature.}

\medskip

We will think that the alternating group $\Alt(N)$ acts
on a set of $N$ points which are arranged into $6$ dimensional
cube of size $K$ and we will identify these points with ordered
$6$-tuples of nonzero elements from the field $\F_{2^{3s}}$.

Let $\rho:\Alt(N) \to U(\mr{H})$ be a fixed unitary representation of
the alternating group and let $v \in \mr{H}$ be
an $\epsilon$-almost invariant unit vector for some generating set
$F$. We will fix the set $F$ and the number $\epsilon$ later.
Without loss of generality, we may assume that
$\mr{H}$ is generated by the orbit of the vector $v$.

\medskip

Let $H_s$ denote the group $\SL_{3s}(\F_2)$. The group $H_s$
has a natural action on the set $V\setminus \{0\}$ of $K$ nonzero
elements of a vector space $V$ of dimension $3s$ over $\F_2$.
The elements of $H_s$ act by even permutations on $V\setminus \{0\}$, because
$H_s$ is a simple group and does not have $\Z/2\Z$ as a factor.
If we identify $V$ with $\F_{2^{3s}}$ then the
existence of a generator for the multiplicative group of $\F_{2^{3s}}$ implies
that some element of $H_s$ acts as a $K$-cycle on $V\setminus\{0\}$.%
\footnote{We can use some other family of groups instead of $\{H_s\}_s$ --- the only
requirements for $\{H_s\}_s$ are that they acting transitively on a set of $K(s)$ points
(where $K(s) \to \infty $ as $s \to \infty$)
and that $\left\{{H_s}^{\times {K(s)^5}}\right\}_s$ can be made a bounded degree expanders.
For example we can use the groups $\SL_3(\F_p)$ acting on ${\F_p}^3 \setminus \{0\}$ or
on the projective plane ${P^2\F_p}$ for different primes $p$-es.}

Let $\Gamma$ be the direct product of $K^5$ copies of the group $H_s$.
The group $\Gamma$ can be embedded into $\Alt(N)$ in $6$ different ways
which we denote by $\pi_i$, $i=1,\dots, 6$. The image of
each copy of $H_s$ under $\pi_i$ acts as $\SL_{3s}(\F_2)$
on a set of $K=2^{3s}-1$
points where all
coordinates but the $i$-th one are fixed. It is clear that
$\Gamma$ contains an abelian subgroup $\overline{\Gamma}$
isomorphic to $(\Z/K\Z)^{\times K^5}$.

Using Theorem~5 from~\cite{KSL3k}, we can find a small generating set
$S$ of the group $H_s$ such that the \KaC\ $\KC(H_s;S) > 1/400$.
This allows us to construct a generating set $\tilde S$ with $40$ elements
of $\Gamma$ with similar properties, i.e., the \KaC\
$$
\KC(\Gamma;\tilde S) > 1/500.
$$

Now we can define the generating set $F_N$ of $\Alt(N)$ such that the
\KaC\ $\KC(\Alt(N);F_N)$ can be estimated --- the set $F_N$ will be the
union of the images of $\tilde S$ under the embeddings $\pi_i$:
$$
F_N = \bigcup_i \pi_i(\tilde S).
$$
The group generated by the set $F_N$ contains the $6$ images
of $\Gamma$ and therefore
is the whole alternating group $\Alt(N)$. From now on, we will assume
that the vector $v\in\mr{H}$ is $\epsilon$-almost invariant with
respect to the set $F_N$.

Using the \KaC\ of $\KC(\Gamma;\tilde S)$, it can be seen that if
$v$ is $\epsilon$-almost invariant vector with respect to the set $F_N$
in some representation $\rho$ of $\Alt(N)$, then $v$ is close to $\Gamma_i$
invariant vector, i.e.,
$$
||\rho(\pi_i(g))v -v || \leq 10^3 \epsilon ||v||,
$$
for all $g \in \Gamma$ and any $i=1,\dots, 6$.
If the diameters of Cayley graphs of $\mr{C}(\Alt(N),\cup \Gamma_i)$ were
bounded independently on $N$, this would gives us that the
representation $\rho$
has an invariant vector, provided that $\epsilon$ is small enough.
Unfortunately this is not the case, because $\Alt(N)$ can not be written
as a product of less than $\log |N|$ abelian subgroups
and each $\Gamma_i$ is a product of less than $20$ abelian subgroups.

Let $\overline{\Gamma_i}$ denote the image of $\overline{\Gamma}$
under the embedding $\pi_i$ and let $C$ be the union of $\overline{\Gamma_i}$.
Thus, 
we may assume that the vector $v$ is almost invariant with respect to the set $C$.

\medskip

As mentioned before we will break the representation $\rho$ into two components.
The space $\mr{H}$ decomposes as a sum of irreducible representations
$$
\mr{H} = \bigoplus_\lambda c_\lambda M^\lambda,
$$
where the sum is over all partitions $\lambda$ of $N$, $M^\lambda$ denotes the
irreducible representations corresponding to the partition $\lambda$
and $c_\lambda$ is the multiplicity of $M^\lambda$ in $\mr{H}$ which is
either $0$ or $1$, since $\mr{H}$ is generated  by 1 element.
Let $\mr{H}_1$ is the sum of all irreducible sub-representations
of $\mr{H}$ which correspond
to the partitions $\lambda=[\lambda_1,\lambda_2,\dots]$ of $N$
with
$$
\lambda_1 < N - K^{5/4},
$$
and $\mr{H}_2$ is the orthogonal complement of $\mr{H}_1$ in $\mr{H}$.

This allows us to decompose the almost invariant vector $v$ as
$v=v_1 + v_2$, where
$v_i \in \mr{H}_i$. We will use two different arguments to
show that the vector $v_1$
is small and that $v_2$ is close to an invariant vector in $\mr{H}_2$.

\medskip

Using the definition of the set $C$, it can be seen that $C^{440}$
acts almost transitively
on the set of all ordered tuples of $K^5/10$ points.
Here $C^{440}$ denotes the set
of all product of less than $440$ elements from the set $C$,
and by almost transitivity we if we are
given $2$ ordered tuples then with large probability $\kappa$
(approaching $1$ as $s\to \infty$), there
exists an element in $C^{440}$ which send one to the other.
Therefore the set $C^{440}$ contains almost all elements in some
conjugacy class $B$ of permutations
in $\Alt(N)$ with at least $K^5/10$ non-fixed points.

The vector $v$ is almost preserved by any element of $C$,
which implies that $v$
is moved a little by most of the elements inside the
conjugacy class $B$, i.e.,
$||\rho(g) v -v|| \leq 440\times 1000 \epsilon $ for any
$g \in B \cap C^{440}$. This gives
$$
\begin{array}{l@{}c@{}l}
\displaystyle
\left|\left| \frac{1}{|B|}\sum_{g\in B} \rho_1(g) v_1 -v_1 \right| \right|
& \,\,{} \leq {}\,\,&
\displaystyle
\frac{1}{|B|}\sum_{g\in B} \left|\left|  \rho(g) v -v \right| \right| \leq \\
& \leq &
\displaystyle \rule[20pt]{0pt}{0pt}
 2(1-\kappa) ||v_1|| +
\frac{\kappa}{|B|}\sum_{g\in B \cap C^{440}}
         \left|\left|  \rho(g) v -v \right| \right| < \\
& < &
\displaystyle
\frac{||v_1||}{6} + 4.5 \times 10^5 \epsilon ,
\end{array}
$$
provided that $\kappa$ is sufficiently close to $1$.
The decomposition of  $\mr{H}_1$ as
$$
\mr{H}_1 = \bigoplus_\lambda c_\lambda M^\lambda,
$$
gives a decomposition of the vector $v_1 = \sum v_\lambda$.
The set $B$ is a conjugacy class therefore
$\frac{1}{|B|}\sum_{g\in B} \rho(g) v_\lambda =
   \bar \chi_{\lambda}(B) v_\lambda $
for any vector $v_\lambda$  in an irreducible representation $M^{\lambda}$.
Here $\bar \chi_{\lambda}(B)$ is the normalized character of the representation
$M^{\lambda}$,
defined as $\bar \chi_{\lambda}(B) := \chi_{\lambda}(B)/ \dim M^{\lambda}$.
Thus we have
$$
\left|\left| \frac{1}{|B|}\sum_{g\in B} \rho_1(g) v_1 \right| \right| \leq
||v_1||\,\, \max_{\lambda}\,\,\,\,| \bar \chi_{\lambda}(B)|,
$$
where the maximums are taken over all partitions which appear in
the representation $\rho_1$, i.e.,
all partitions $\lambda$ with  $\lambda_1 < N - K^{5/4}$.
There are various estimates of value of the normalized characters of the
symmetric/alternating groups. Applying the bounds from~\cite{Ro} Theorem~1 yields
$$
\begin{array}{l@{}c@{}l}
\displaystyle
\max_{\lambda} |\bar \chi_{\lambda}(B)|
& {} \leq {}&
\displaystyle
\max_{\lambda} \left\{ \max
        \left\{ \lambda_1/N ,q \right\}^{c\, \mbox{\tiny supp}|B|} \right\}
\leq \\
& \leq &
\displaystyle \rule[20pt]{0pt}{0pt}
\left(1-\frac{K^{\frac{5}{4}}}{K^6}\right)^{\frac{c K^5}{10} } \leq
  \exp\left(-\frac{cK^{\frac{1}{4}}}{10}\right)
  < \frac{1}{3},
\end{array}
$$
where $c$ and $q$ are universal constants. The last inequality is valid only if $K$ is large enough.
The two inequalities above imply that
\begin{equation}
\label{v1}
||v_1|| < 9 \times 10^{5} \epsilon.
\end{equation}

\medskip

The above argument does not work for the representation $\mr{H}_2$, because the first
part of the partition $\lambda$ can be close to $N$. This means that
the sum
$$
\lambda_2+ \lambda_3 + \dots < K^{5/4}
$$
is small. Thus
$\mr{H}_2$ can be embedded in the representation $M$ arising from
the action of $\Alt(N)$ on set of all ordered tuples of size $K^{5/4}$.
Let $E$ denote the set of ordered tuples of size $K^{5/4}$
which is a basis of $M$.

We have $K^{K^5} \gg N^{K^{5/4}}$, i.e.,
number of elements in the set $C$ is much larger than size of the set $E$.
Using this inequality and the definition of the set $S$,
it can be shown that the random walk on set $E$,
where the moves are given by the permutations from some subset
$\widetilde{C} \subset C^6$, mixes in a few steps
independent of $N$.
Therefore if we define the operator $\Delta$ on $M$ defined by
$$
\Delta:=\frac{1}{|\widetilde{C}|}\sum_{g\in \widetilde{C}} \rho_2(g)
$$
then $\Delta^8$ has a single eigenvalue $1$ with eigenvectors the
invariant vectors in $M$, and
all other eigenvalues are less than $1/2$ by absolute value.
Thus we have:
$$
\left|\left| \Delta^8 v_2 -v _2 \right|\right| \geq
\frac{1}{2}||v_2 - v_{||}||,
$$
where $v_{||}$ be the projection of $v_2$ onto the space of
all invariant vectors in $M$.
On the other hand
$$
\left|\left| \Delta^8 v_2 -v _2 \right|\right| \leq
 8 \left|\left| \Delta v_2 -v _2 \right|\right| \leq
\frac{8}{|\widetilde{C}|}\sum_{g\in \widetilde{C}}
         || \rho_2(g) v_2 - v_2 || \leq 48 \times 1000 \epsilon,
$$
which gives that
\begin{equation}
\label{v2}
||v_2 - v_{||}|| < 10^5 \epsilon.
\end{equation}

The inequalities (\ref{v1}) and (\ref{v2})
imply that
$$
||v-v_{||}||  \leq  || v_1 || + || v_2 - v_{||}|| < 10^6 \epsilon.
$$
In particular, if $\epsilon$ is small enough
then the vector $v_{||}$ is not zero, which show that there
exists invariant vectors in the representation $\mr{H}$.
Thus, we have shown that
$$
\KC(\Alt(N);F_N) \geq 10^{-6},
$$
which concludes the proof of Theorem~\ref{main}.
$\square$

\medskip

\noindent
\textbf{Proof of Theorem~\ref{main1}:}
By Theorem~\ref{main} the alternating groups $\Alt(n_s)$
are expanders with respect
to some generating set $F_{n_s}$ for $n_s=\left(2^{3s}-1\right)^6$.
The sequence $\left\{n_s\right\}_{s=50}^{\infty}$
grows exponentially, therefore for
any sufficiently large $n$ there exists $s$
such that $1< n/n_s < 10^6$. The group $\Alt(n)$ can be written as a product of
fixed number (less then $~10^{9}$) of copies of $\Alt(n_s)$
embedded in $\Alt(n)$. Using the images of
the sets $F_{n_s}$ one can construct a generating set $F_n$ such that
$$
\KC\left(\Alt(n); F_n\right) \geq 10^{-15}
$$
and $|F_n| \leq 10^{10}$, which completes the proof of Theorem~\ref{main1},
the construction of the
generating sets $\tilde F_n$ of the symmetric groups is similar. $\square$

\medskip

The generating set $S$ of $\SL_{3s}(\F_2)$ can be defined so that the elements
of $S$ are involutions. This allows us to construct an expanding generating
set $F_n$ consisting only of involutions.

The bounds for the size of the generating set $F_n$ and the \KaC\ in the proof Theorem~\ref{main1}
can be significantly improved ---
it is possible to construct a $10$ element generating set
$\bar F_n$  consisting of
involutions such that
$\KC\left(\Alt(n); \bar F_n\right) \geq 10^{-8}$, provided that
$n$ is large enough. N. Nikolov~\cite{Np} suggested that it is possible
to further improve the bound for \KaC\ by
a factor of two by using the groups $\SL_3(\F_q)$ instead of $\SL_{3s}(\F_2)$,
but this will increase the size of the generating set.

\bigskip

\noindent
\textbf{Acknowledgements: }
I wish to thank A.~Lubotzky and N.~Nikolov 
for their encouragement and the useful discussions during the work on this project.
I am grateful to Y.~Shalom and E.~Zelmanov for introducing me to the subject.


\texttt{\\Martin Kassabov, \\
Cornell University, Ithaca, NY 14853-4201, USA. \\
\emph{e-mail:} kassabov@math.cornell.edu }
\end{document}